\documentclass{commat}

%%% AUTHOR'S PACKAGES %%%
\usepackage{graphicx}
\usepackage{float}

\title{%
    Qualitative analysis of strictly non-Volterra quadratic
dynamical systems with continuous time
    }

\author{%
    Xaydar Raupovich Rasulov
    }

\affiliation{
    \address{Rasulov Xaydar Raupovich, Department of Mathematical analysis, Bukhara state university, Institute
of Mathematics, Bukhara branch, Uzbekistan
Academy of Sciences, Bukhara, Uzbekistan
        }
    \email{%
    email@xrasulov71@mail.ru
    }
    }

\abstract{%
    In this article, a continuous analogue of strictly non-Volterra quadratic dynamical systems with continuous time and points of equilibrium is investigated, a phase portrait of the system is constructed, numerical solutions are found, and a comparative analysis is carried out with a particular solution of the system.
    }



\VOLUME{30}
\NUMBER{1}
\firstpage{239}
\DOI{https://doi.org/10.46298/cm.10528}

\begin{paper}

\section{Introduction}

The main theoretical and practical studies on the study of the dynamics of a free population, various problems of physics and economics are reduced to dynamic systems. In this, a special role is played by the quadratic stochastic operators. In this regard, quadratic operators attract the attention of specialists in various fields of mathematics and its applications (see for example \cite{Bern2}, \cite{Ulam21}). Therefore, studies of quadratical stochastic operators remain relevant.

The concept of a quadratic stochastic operator was first formulated in the article \cite{Bern2}. S. Ulam \cite{Ulam21} posed the problem of studying the behavior of the trajectories of quadratic stochastic operators. This problem is mainly solved for Volterra operators (see \cite{Gani5}, \cite{Gani6}, \cite{Gani7}, \cite{Gani8}, \cite{Mukh14}, \cite{Gani4}) of discrete time.
In \cite{Tian19}, a definition of strictly non-Volterra quadratic stochastic operators is introduced, which is a subclass of non-Volterra operators. But the class of non-Volterra operators has been little studied.

Depending on the problem, either continuous time can be considered (when the states of the system are of interest at each moment), or discrete (when the states of the system are of interest at separate isolated moments of time). System behavior depends on parameters, the exact values of which are often unknown. It is fundamentally impossible to solve such problems numerically for all possible parameter values. Therefore, methods are needed that will make it possible to analyze the behavior of solutions to a dynamical system without using a computer.

As we noted above, quadratic stochastic operators of the Volterra class in discrete time have been studied quite deeply. Here are some of them \cite{Rozi17}, \cite{Gani9}, \cite{Kest12}, \cite{Lyub13}, \cite{Rozi18}. Also, quadratic stochastic operators with continuous time, that is, various modifications of the predator and prey model and other models, have been studied by many authors in \cite{Ange1}, \cite{Tian19} and \cite{Huan10}.

In turn, the class of strictly non-Volterra quadratic stochastic operators with continuous time (hereinafter referred to as quadratic dynamical systems) has been studied relatively less. Since there is no general theory studying such operators. Note the works \cite{Rasu16}.

The most useful aspect of the qualitative theory of dynamical systems, be more correct, dynamical systems with continuous time, is that many important properties of solutions can be predicted in advance without having explicit solutions of the equations.

In this article, we study a qualitative analysis of a quadratic dynamical system, a continuous analogue of the non quadratic stochastic Volterra operator from \cite{Mukh15}. The system is shown to come to a simple differential equation, solutions are found, and they are studied separately. Equilibrium points of the system are found and a phase portrait is given.

Analytical solutions of the basic system in some assumptions have been found. In contrast to the results of the study of discrete time dynamic systems, the trajectory of the solutions is a clearly expresses of the aspiration to the equilibrium point at an exponential speed. With the help of computer calculations, numerical solutions were compared with analytical solutions, graphs were drawn, and trajectories were found to overlap over time.

\section{Formulation of the problem}

In this article, we study a continuous analogue of one quadratic stochastic operator from \cite{Mukh15}, which in our case has the form:
\begin{equation}\label{eq1}
\begin{cases}
\dot{x}_0=\frac{1}{2}x_1^2+\frac{1}{2}x_2^2+2x_1 x_2-x_0=f_0(x_0,x_1,x_2 ),\\
\dot{x}_1=\frac{1}{2}x_0^2+\frac{1}{2}x_2^2+2x_0 x_2-x_1=f_1(x_0,x_1,x_2 ),\\
\dot{x}_2=\frac{1}{2}x_0^2+\frac{1}{2}x_1^2+2x_0 x_1-x_2=f_2(x_0,x_1,x_2 )
\end{cases}
\end{equation}
or in vector form $\dot{x}(t)=f(x(t) )$, where $x(t)=(x_0 (t), x_1 (t), x_2 (t) ),$ is the state of some system at the instant of continuous time at $t\geq0,  x_0 (t)\geq0,  x_1 (t)\geq0, x_2 (t)\geq0$ and $x_0 (t)+x_1 (t)+x_2 (t)=1.$

\section{Main results}

Adding all the equations of system \eqref{eq1} and denoting $x_0+x_1+x_2=X$ we obtain an ordinary differential equation
\begin{equation}
X^{'}=X^2-X,\label{eq2}
\end{equation}
which is a type with shared variables.
	
Equation \eqref{eq2} has the following solution:
\begin{equation}
X= \frac{1}{1-Ce^t},\label{eq3}
\end{equation}
where $C=const\ne 0.$
For a fixed negative $C$, formula \eqref{eq3} gives one solution located in the line $0<X<1$. For a fixed positive $C$, formula \eqref{eq3} defines two solutions, one of which, defined on the interval $-\infty<t<-lnC$, is located in the half-plane $X>1$, and the other, defined on the interval $-lnC<t<+\infty,$ is located in the half-plane $X<0.$ In addition to solutions located in these regions indicated above, equation \eqref{eq2} has two more solutions $X\equiv0$ and $X\equiv1$, which are formally obtained from \eqref{eq3} with $C=0$ and $C=\infty$.

This means that the phase portrait of equation \eqref{eq2} consists of five phase curves: two singular points $0$ and $1$, two rays $(-\infty,0)$ and $(1,+\infty)$, and an interval $(0,1)$. Here, the point $X=0$ is stable, and the point $X=1$ is an unstable equilibrium.

Let's take a look at each solution separately. Let us show that only one solution $X=1$ corresponds to the considered system \eqref{eq1}.

So, $X=0$ i.e., $x_0+x_1+x_2=0$. From here, we find $x_0=-x_1-x_2$ and put on the second and third equations of system \eqref{eq2}
\begin{equation}
\begin{cases}
\dot{x}_0=\frac{1}{2}x_1^2+\frac{1}{2}x_2^2+2x_1 x_2-x_0 ,\\
\dot{x}_1=\frac{1}{2}x_1^2-x_2^2-x_1 x_2-x_1 ,\\
\dot{x}_2=-x_1^2+\frac{1}{2}x_2^2+2x_1 x_2-x_2 .
\end{cases}\label{eq4}
\end{equation}

After some transformations of the second and third equations of system \eqref{eq4}, we obtain an ordinary Abel differential equation of the second kind with respect to $x_1,$ which does not admit a solution in quadratures \cite{Kamk11},\cite{Tokm20}, \cite{Zait22}. Hence we can say that system \eqref{eq4} in the general case has no solution in the analytical form.

\begin{definition}
Equilibrium positions of system \eqref{eq1} are such points $x^{*}(t)$ of the phase space, which $f(x^{*}(t))=0$. Obviously, $x^{*}(t)$  is a solution of system \eqref{eq1}, since $\dot{x}^{*}=0$ .
\end{definition}

System \eqref{eq4} has 4 equilibrium points: $N_{1}(0,0,0)$, $N_{2}(\frac{4}{3}, -\frac{2}{3},-\frac{2}{3})$, $N_{3}(-\frac{2}{3},-\frac{2}{3},\frac{4}{3})$ and $N_{4}(-\frac{2}{3},\frac{4}{3},-\frac{2}{3})$, of which only $N_{1}$
is located on the border of the first octant. If $x_{0}(t)\ge 0$, $x_{1}(t)\ge 0$, $x_{2}(t)\ge 0$ and $x_{0}(t)+x_{1}(t)+x_{2}(t)=1$, for all $t\ge 0$, then the solution $X=0$ cannot describe the original problem.

Although, the solution of the system located in the line $0<X<1$
$$X=\frac{1}{1-Ce^{t}}$$
for a fixed negative C satisfies the condition for $t\ge 0$, $
x_{0}(t)\ge 0$, $x_{1}(t)\ge 0$, $x_{2}(t)\ge 0$, but does not satisfy the
condition $x_{0}(t)+x_{1}(t)+x_{2}(t)=1$. And the solution of the
system located on the lines $X<0$ and $X>1$ does not satisfy any
condition which, for $t\ge 0$, $x_{0}(t)\ge 0$, $x_{1}(t)\ge 0$, $x_{2}(t)\ge
0$ and $x_{0}(t)+x_{1}(t)+x_{2}(t)=1$.

We can conclude that only one solution $X\equiv 1$ of equation \eqref{eq2}
describes the considered problem \eqref{eq1}.

From the condition $x_{0}(t)+x_{1}(t)+x_{2}(t)=1$ we find $
x_{2}=1-x_{0}-x_{1}$ and assuming the first and second equations of
system \eqref{eq1} we get:
\begin{equation}
\begin{cases}
\dot{x}_{0}=\frac{1}{2}x_{0}^{2}-x_{1}^{2}-x_{0}x_{1}-2x_{0}+x_{1}+\frac{1}{2},\\
\dot{x}_{1}=-x_{0}^{2}+\frac{1}{2}x_{1}^{2}-x_{0}x_{1}+x_{0}-2x_{1}+\frac{1}{2},\\
\dot{x}_{2}=\frac{1}{2}x_{0}^{2}+\frac{1}{2}x_{1}^{2}+2x_{0}x_{1}-x_{2}.
\end{cases}\label{eq5}
\end{equation}

By substitution, the first two equations of the system \eqref{eq5}
$$
\begin{cases}
3x_{0}+3x_{1}=v_{0}, \\
x_{0}-x_{1}=v_{1},
\end{cases}
$$
are reduced to the system of equations:
$$
\begin{cases}
v_{0}^{'}=\frac{3}{4}v_{1}^{2}-\frac{1}{4}v_{0}^{2}-v_{0}+3, \\
v_{1}^{'}=\frac{1}{2}v_{0}v_{1}-3v_{1}.
\end{cases}
$$

Finding $v_{0}$  and, accordingly, $v_{0}^{'}$ from the second
equation of this system and substituting into the first equation and
setting $v_{1}^{'}=Y, v_{1}^{''}=YY^{'},$ we obtain the Abel equation
of the second type:
$$
v_{1}YY^{'}=\frac{1}{2}Y^{2}-4v_{1}Y+\frac{3}{8}v_{1}^{4}-6v_{1}^{2},
$$
which does not admit a solution by quadratures \cite{Kamk11},\cite{Tokm20}, \cite{Kest12}. Hence we can say that system \eqref{eq1} in the general case has no
solution in the analytical form. Below, we show that system \eqref{eq1}, under
some assumptions about $x_{0}(t)$, $x_{1}(t)$, $x_{2}(t)$ admits solutions
in quadratures.

 We find the equilibrium positions of system \eqref{eq1}, that is,
consider the equation $f(x)=0.$
$$
\begin{cases}
\frac{1}{2}x_{1}^{2}+\frac{1}{2}x_{2}^{2}+2x_{1}x_{2}-x_{0}=0, \\
\frac{1}{2}x_{0}^{2}+\frac{1}{2}x_{2}^{2}+2x_{0}x_{2}-x_{1}=0, \\
\frac{1}{2}x_{0}^{2}+\frac{1}{2}x_{1}^{2}+2x_{0}x_{1}-x_{2}=0.
\end{cases}
$$

\hspace*{-2pt}This  system  has  8  points  of  equilibrium  position:  $M_{1}(0,0,0)$,  $M_{2}(\frac{1}{3}, \frac{1}{3}, \frac{1}{3})$, $M_{3}(-1,-1,3)$,
$M_{4}(-\frac{2}{3},-\frac{2}{3},\frac{4}{3})$,   $M_{5}(-1,3,-1)$,   
 $M_{6}(3,-1,-1)$, $M_{7}(-\frac{2}{3},\frac{4}{3},-\frac{2}{3})$ and 
 $M_{8}(\frac{4}{3},-\frac{2}{3},-\frac{2}{3})$, where c is an
arbitrary constant.

Considering the conditions of the problem $x_{0}(t)\ge 0$, $x_{1}(t)\ge 0$, $x_{2}(t)\ge 0$ and $x_{0}(t)+x_{1}(t)+x_{2}(t)=1$,
then we obtain a single point of equilibrium position $M_{2}(\frac{1}{3},\frac{1}{3},\frac{1}{3})$, corresponding to system
\eqref{eq1}.

Let us investigate the stability of the solution at the point $
M_{2}(\frac{1}{3},\frac{1}{3},\frac{1}{3})$ of system \eqref{eq1}. We
calculate at the point $M_{2}$ the Jacobian:
$$A=\begin{pmatrix}
\frac{\partial f_{0}}{\partial x_{0}} & \frac{\partial f_{0}}{\partial
x_{1}} & \frac{\partial f_{0}}{\partial x_{2}}\\
\frac{\partial f_{1}}{\partial x_{0}} & \frac{\partial f_{1}}{\partial
x_{1}} & \frac{\partial f_{1}}{\partial x_{2}}\\
\frac{\partial f_{2}}{\partial x_{0}} & \frac{\partial f_{2}}{\partial
x_{1}} & \frac{\partial f_{2}}{\partial x_{2}}
\end{pmatrix}_{M_{2}(\frac{1}{3}, \frac{1}{3}, \frac{1}{3})}=\begin{pmatrix}
-1 & 1 & 1 \\
1 & -1 & 1 \\
1 & 1 & -1
\end{pmatrix},
$$
where $det\vert A\vert =4\ne 0$.
Linearized system \eqref{eq1} has the following form:
\begin{equation}
\begin{cases}
x_{0}^{'}=-x_{0}+x_{1}+x_{2}-\frac{1}{3}, \\
x_{1}^{'}=x_{0}-x_{1}+x_{2}-\frac{1}{3}, \\
x_{2}^{'}=x_{0}+x_{1}-x_{2}-\frac{1}{3}.
\end{cases}\label{eq6}
\end{equation}

Assuming $u=x_{0}-\frac{1}{3}$, $v=x_{1}-\frac{1}{3}$, $w=x_{2}-\frac{1}{3}
$, we get
\begin{equation}
A=\begin{cases}
u^{'}=-u+v+w, \\
v^{'}=u-v+w, \\
w^{'}=u+v-w.
\end{cases}\label{eq7}
\end{equation}

Here, the corresponding matrix is
$\begin{pmatrix}
-1 & 1 & 1 \\
1 & -1 & 1 \\
1 & 1 & -1
\end{pmatrix}
$, its determinant is $det A = 4 \ne 0$, and its eigenvalues are $\lambda_{1, 2}=-2$, $\lambda_{3}=1$. This means that the point of equilibrium
position $M_{2}$ is an unstable rest point.

Let $n_{+}, n_{0}, n_{-}-$ be the number of eigenvalues of $A
$ (taking into account their multiplicity) with positive, equal to zero
and negative real parts, respectively (see \cite{Brat3}).

\begin{definition} The equilibrium position of the dynamical
system \eqref{eq7} is called hyperbolic if $n_{0}=0$, that is, there are no
eigenvalues located on the imaginary axis. A hyperbolic equilibrium is
called a hyperbolic saddle if $n_{+}n_{-}\ne 0.$ Hence, on the other
hand, the equilibrium position $M_{2}$ is a hyperbolic saddle.
The general solution of system \eqref{eq6} has the form:
\begin{equation}
\begin{cases}
x_{0}=\frac{1}{3}-(C_{0}+C_{1})e^{-2t}+C_{2}e^{t}, \\
x_{1}=\frac{1}{3}+C_{0}e^{-2t}+C_{2}e^{t}, \\
x_{2}=\frac{1}{3}+C_{1}e^{-2t}+C_{2}e^{t},
\end{cases}\label{eq8}
\end{equation}
where $C_{0},C_{1},C_{2}=const$.
\end{definition}

If $x_{0}(t)\ge 0$, $x_{1}(t)\ge 0$, $x_{2}(t)\ge 0$ and $x_{0}(t)+x_{1}(t)+x_{2}(t)=1$, then we have to put $C_{2}=0$ and the solution will look like:
\begin{equation}
\begin{cases}
x_{0}=\frac{1}{3}-(C_{0}+C_{1})e^{-2t}, \\
x_{1}=\frac{1}{3}+C_{0}e^{-2t}, \\
x_{2}=\frac{1}{3}+C_{1}e^{-2t}.
\end{cases}\label{eq9}
\end{equation}

It follows that, as $t\to +\infty $ the solution of system \eqref{eq5} $
(x_{0}(t),x_{1}(t),x_{2}(t))$ tends to the point $
M_{2}(\frac{1}{3},\frac{1}{3},\frac{1}{3})$ exponentially fast.

There are two constants in \eqref{eq9} that cannot uniquely determine
the solution to the Cauchy problem for system \eqref{eq6}. If we put $
x_{1}(t)=x_{2}(t),$ then the Cauchy problem for system \eqref{eq6} has a
unique solution.

In studying the Cauchy problem for equation \eqref{eq1}, proceeding
from the above, we assume $x_{1}(t)=x_{2}(t)$. Substituting $
x_{0}=1-x_{1}-x_{2}$ into the second and third equations of system
\eqref{eq1}, we obtain:
\begin{equation}
\begin{cases}
\dot{x}_{0}=\frac{1}{2}x_{1}^{2}+\frac{1}{2}x_{2}^{2}+2x_{1}x_{2}-x_{0},\\
\dot{x}_{1}=\frac{1}{2}x_{1}^{2}-x_{2}^{2}-x_{1}x_{2}-2x_{1}+x_{2}+\frac{1}{2},\\
\dot{x}_{2}=-x_{1}^{2}+\frac{1}{2}x_{2}^{2}-x_{1}x_{2}+x_{1}-2x_{2}+\frac{1}{2}.
\end{cases}\label{eq10}
\end{equation}

Assume $x_{1}=x_{2}=v$. Then the second and third equations
take the form:
\begin{equation}
v^{'}=-\frac{3}{2}v^{2}-v+\frac{1}{2}.
\label{eq11}
\end{equation}

The phase portrait of equation \eqref{eq11} consists of five phase curves: two
singular points $v=-1$ and $v=\frac{1}{3}$, two rays $(-\infty
,-1)$ and $(\frac{1}{3},+\infty )$ and an interval $(-1,\frac{1}{3})
$. Here, the point $v=-1$ is stable, and the point $v=\frac{1}{3}$
is an unstable equilibrium.

We solve equation \eqref{eq11} and put the solution in the first
equation of system \eqref{eq10}. Considering $x_{0}(t)\ge 0$, $x_{1}(t)\ge 0$, $x_{2}(t)\ge 0$ and $x_{0}(t)+x_{1}(t)+x_{2}(t)=1$, we find the
general solution of the system \eqref{eq10}.
The solution of system \eqref{eq10} in the line $0<x_{1}(t)$,
$x_{2}(t)<\frac{1}{3}$ has the following form:
\begin{equation}
\begin{cases}
x_{0}=\frac{1}{3}+\frac{8Ce^{-2t}}{3(3+Ce^{-2t})}+C_{2}e^{-t}, \\
x_{1}=\frac{1-Ce^{-2t}}{3+Ce^{-2t}}, \\
x_{2}=\frac{1-Ce^{-2t}}{3+Ce^{-2t}},
\end{cases}\label{eq12}
\end{equation}
and in the line $\frac{1}{3}<x_{1}(t)$, $x_{2}(t)<1$
\begin{equation}
\begin{cases}
x_{0}=\frac{1}{3}-\frac{8Ce^{-2t}}{3(3-Ce^{-2t})}+C_{3}e^{-t}, \\
x_{1}=\frac{1+Ce^{-2t}}{3-Ce^{-2t}}, \\
x_{2}=\frac{1+Ce^{-2t}}{3-Ce^{-2t}},
\end{cases}\label{eq13}
\end{equation}
where $C, C_{2}, C_{3}=const.$

From \eqref{eq12} and \eqref{eq13} we obtain, respectively, that $
x_{0}(t)+x_{1}(t)+x_{2}(t)=1+C_{2}e^{-t}$ and $
x_{0}(t)+x_{1}(t)+x_{2}(t)=1+C_{3}e^{-t}$. Hence it
follows that the sum $x_{0}(t)+x_{1}(t)+x_{2}(t)$ tends to unity
exponentially fast as $t\to +\infty $, and not separately (see
\cite{Mukh15}).
Now, consider the case when $x_{0}(t)=x_{1}(t)=x_{2}(t)=X$.
Substituting this in \eqref{eq2}, taking into account the invariance of the
equations of the system, we obtain the ordinary differential equation
\begin{equation}
X^{'}=3X^{2}-X,\label{eq14}
\end{equation}
which is of the shared variable type too.
As above, it can be shown that equation \eqref{eq14} has the following
solution:
\begin{equation}
X=\frac{1}{3-Ce^{t}}, \label{eq15}
\end{equation}
where $C=const\ne 0.$
For a fixed negative $C$, formula \eqref{eq15} gives one solution located
in the line $0<X<\frac{1}{3}$. For a fixed positive $C$, formula
\eqref{eq15} defines two solutions, one of which, defined on the interval $
-\infty <t<\ln \frac{3}{C}$, is located in the half-plane $
X>\frac{1}{3}$, and the other, defined on the interval $\ln
\frac{3}{C}<t<+\infty $, is located in the half-plane $X<0$. In
addition to solutions located in these regions indicated above, equation
\eqref{eq14} has two more solutions $X\equiv 0$ and $X\equiv \frac{1}{3},$
which are formally obtained from \eqref{eq15} with $C=0$ and $C=\infty.$

The phase portrait of equation \eqref{eq14} consists of five phase curves: two
singular points 0 and $\frac{1}{3}$, two rays $(-\infty ,0)$ and $
(\frac{1}{3},+\infty )$ and an interval $(0,\frac{1}{3})$. Here, the
point $X=0$ is stable, and the point $X=\frac{1}{3}$ is an unstable
equilibrium.

Taking into account, $x_{0}(t)\ge 0,\quad x_{1}(t)\ge 0,\quad x_{2}(t)\ge 0$ and
$x_{0}(t)+x_{1}(t)+x_{2}(t)=1$, as above, we obtain that in this case
$x_{0}(t)=x_{1}(t)=x_{2}(t)=\frac{1}{3},$ i.e. the considered
dynamical system is at rest.

As mentioned above, we failed to find an analytical solution to
the Cauchy problem for the system \eqref{eq1} in the general case. In this
regard, with the help of the MathCAD program, solutions of the Cauchy
problem for system \eqref{eq1} were found and the phase portrait of the
trajectory (with an accuracy of 0.001) was compiled in fig. 1,2.
The phase portrait of system \eqref{eq1} is as follows:

\begin{figure}[H]
\begin{center}
\includegraphics[width=0.5\textwidth]{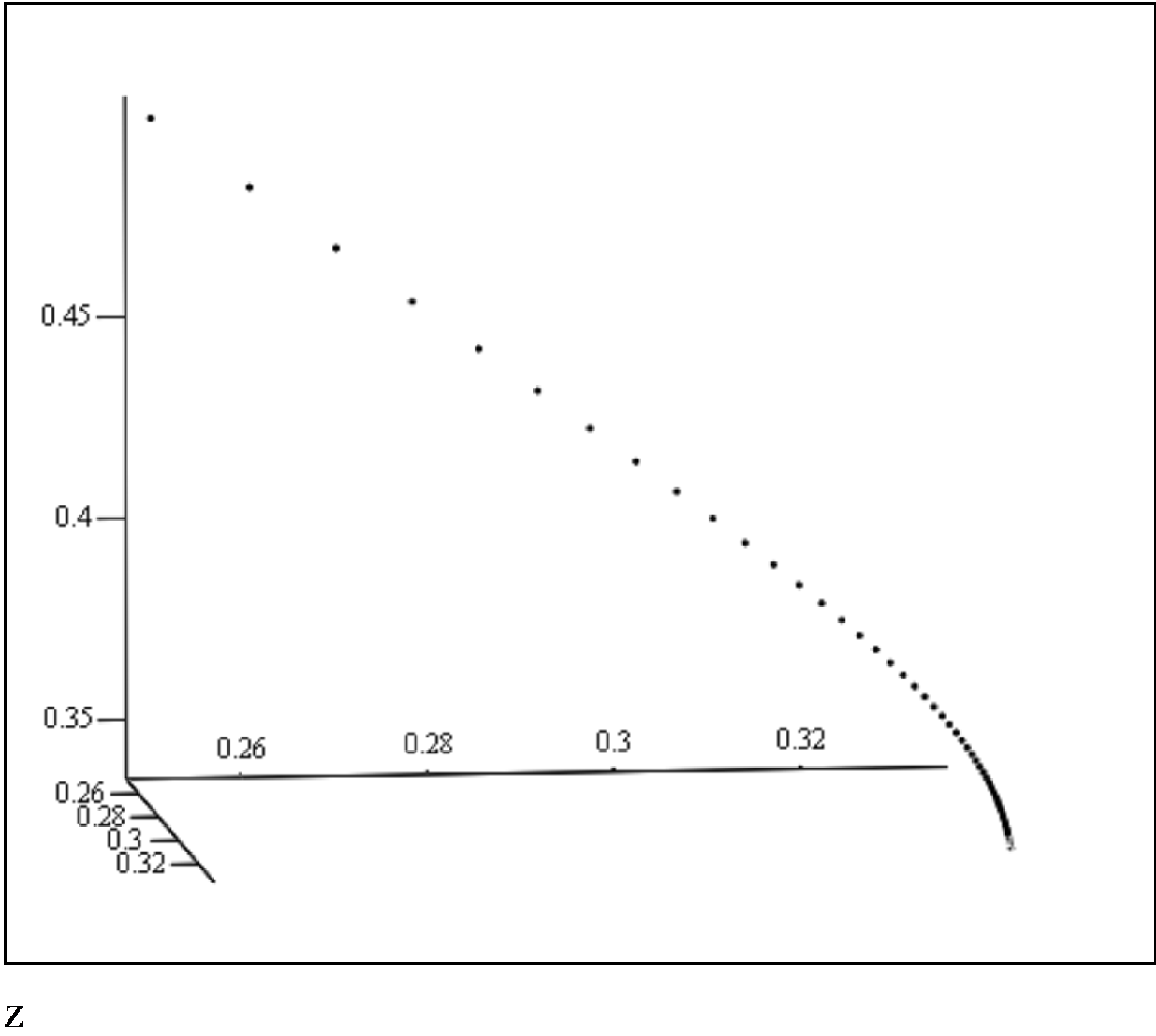}%
\caption{Initial conditions $x_{0}(0)=0.5, x_{1}(0)=0.25,
x_{2}(0)=0.25,$ $Z=(x_{2},x_{1},x_{0}$) (in the line $ 0<x_{1}(t)$, $x_{2}(t)<\frac{1}{3}$, $x_{1}(t)=x_{2}(t))$.}
\end{center}
\end{figure}

\begin{figure}[H]
\begin{center}
\includegraphics[width=0.5\textwidth]{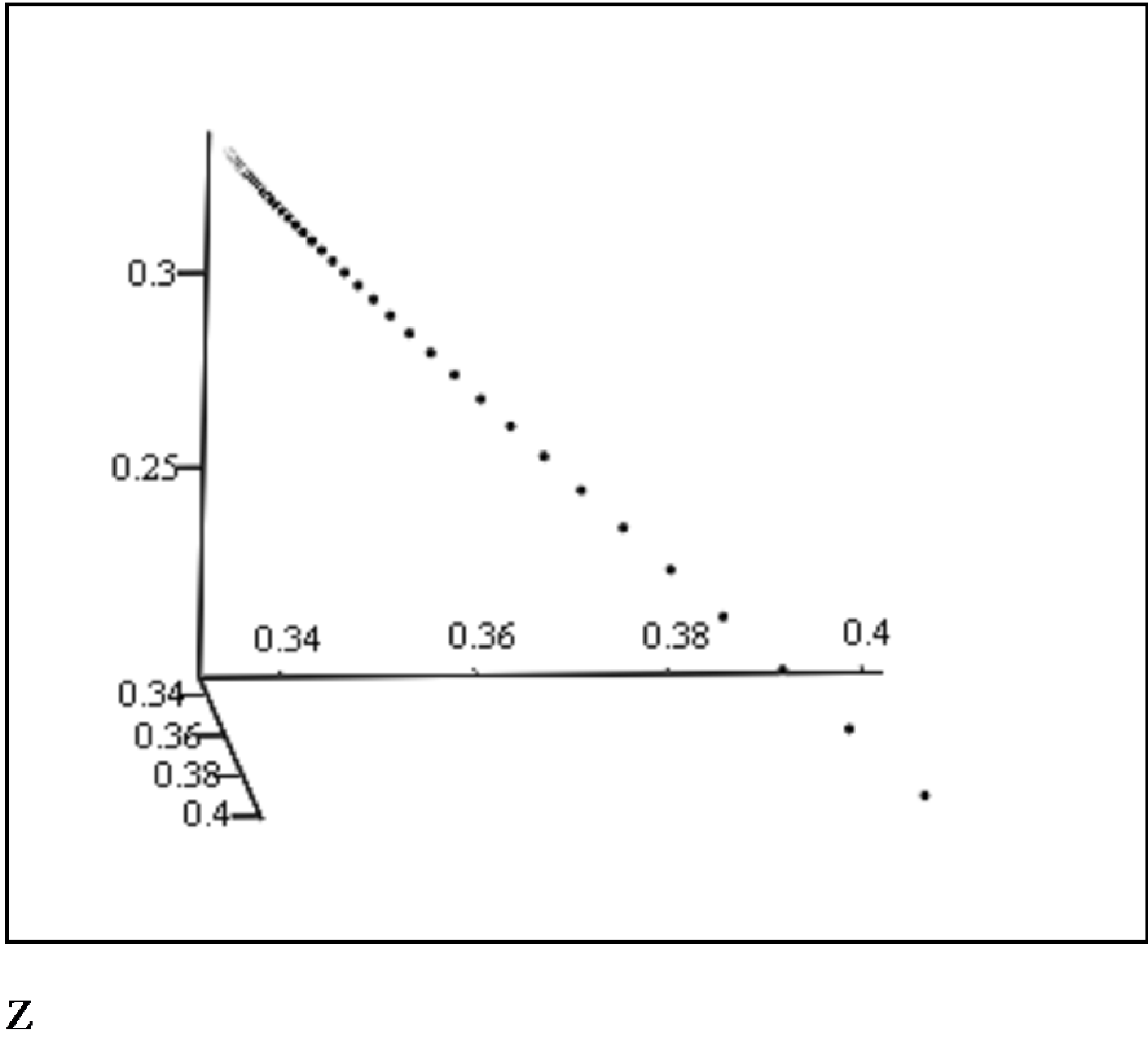}%
\caption{Initial conditions $x_{0}(0)=0.2,\quad x_{1}(0)=0.4$, $x_{2}(0)=0.4$,$Z=(x_{2}, x_{1}, x_{0}$) (in the line $ \frac{1}{3}<x_{1}(t)$, $x_{2}(t)<1$, $x_{1}(t)=x_{2}(t))$.}
\end{center}
\end{figure}

Also, for cases where $x_{1}(t)=x_{2}(t)$ and $0<x_{1}(t)$, $x_{2}(t)<\frac{1}{3}$ and $\frac{1}{3}<x_{1}(t), x_{2}(t)<1$ compared
solutions of system \eqref{eq1}, calculated using the MathCAD program with
\eqref{eq12} and \eqref{eq13}, respectively (fig. 3-6).

To compare the numerical and analytical solutions \eqref{eq12} and \eqref{eq13} on
the graph, we denote by $x_{0}(t), x_{1}(t),x_{2}(t)$ numerical
solutions, and by $y(t), u(t)$ and $w(t)$ analytic solutions from
\eqref{eq12} and \eqref{eq13}, respectively.

As a result of the research, it was found that the difference between
the numerical solutions \eqref{eq1} and \eqref{eq12}, \eqref{eq13} does not exceed 0.001,
and for $t\ge 4$, the numerical and analytical solutions almost
coincide and together tend to the rest point $ M_{2}(\frac{1}{3}, \frac{1}{3}, \frac{1}{3})$.
\begin{figure}[H]
\begin{center}
\includegraphics[width=0.5\textwidth]{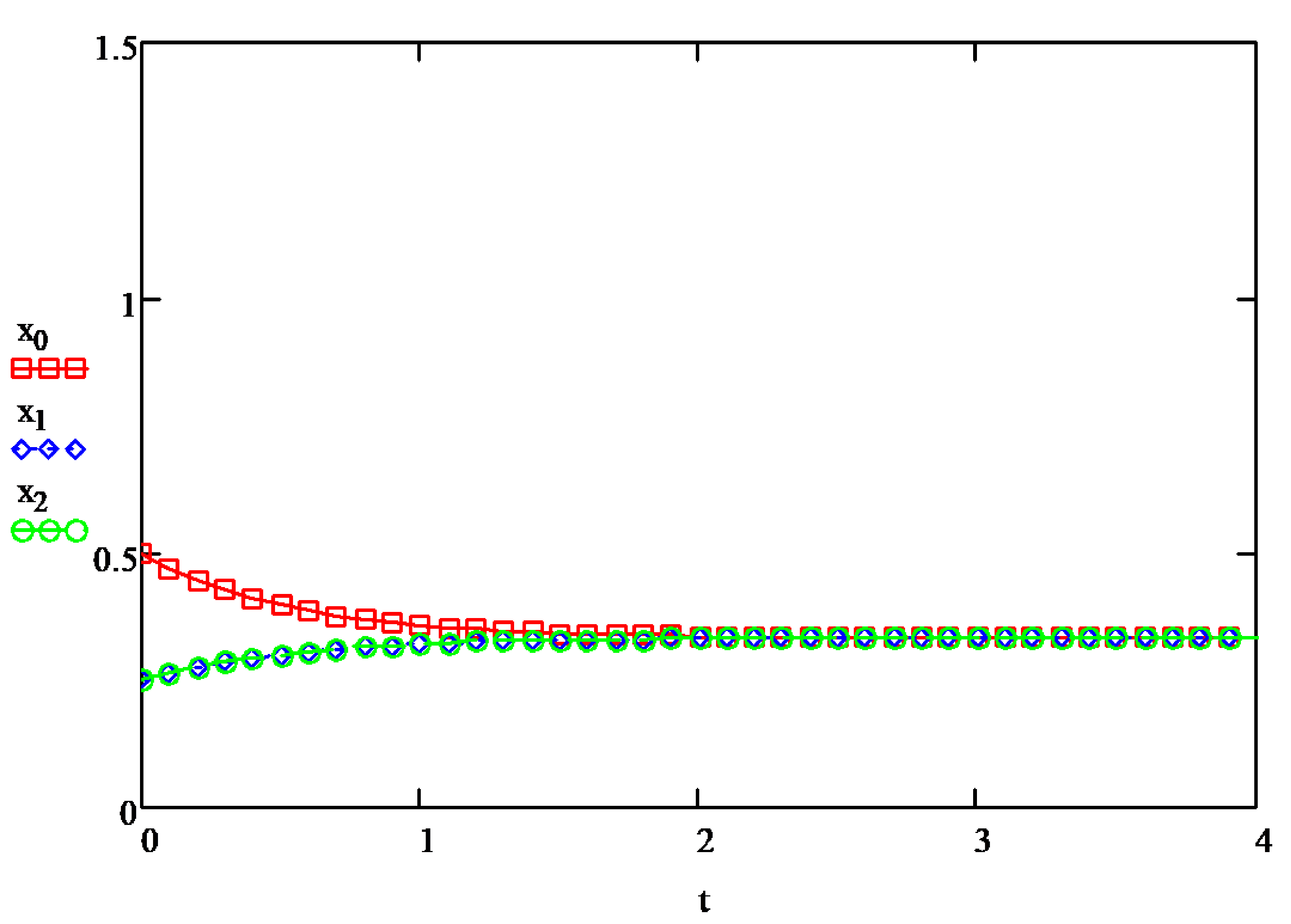}%
\caption{The graph of the numerical solution of system \eqref{eq1} with the
initial conditions $x_{0}(0)=0.5$, $x_{1}(0)=0.25$, $x_{2}(0)=0.25$ (in
the line $0<x_{1}(t)$, $x_{2}(t)<\frac{1}{3}$, $x_{1}(t)=x_{2}(t))$,}
\end{center}
\end{figure}

\begin{figure}[H]
\begin{center}
\includegraphics[width=0.5\textwidth]{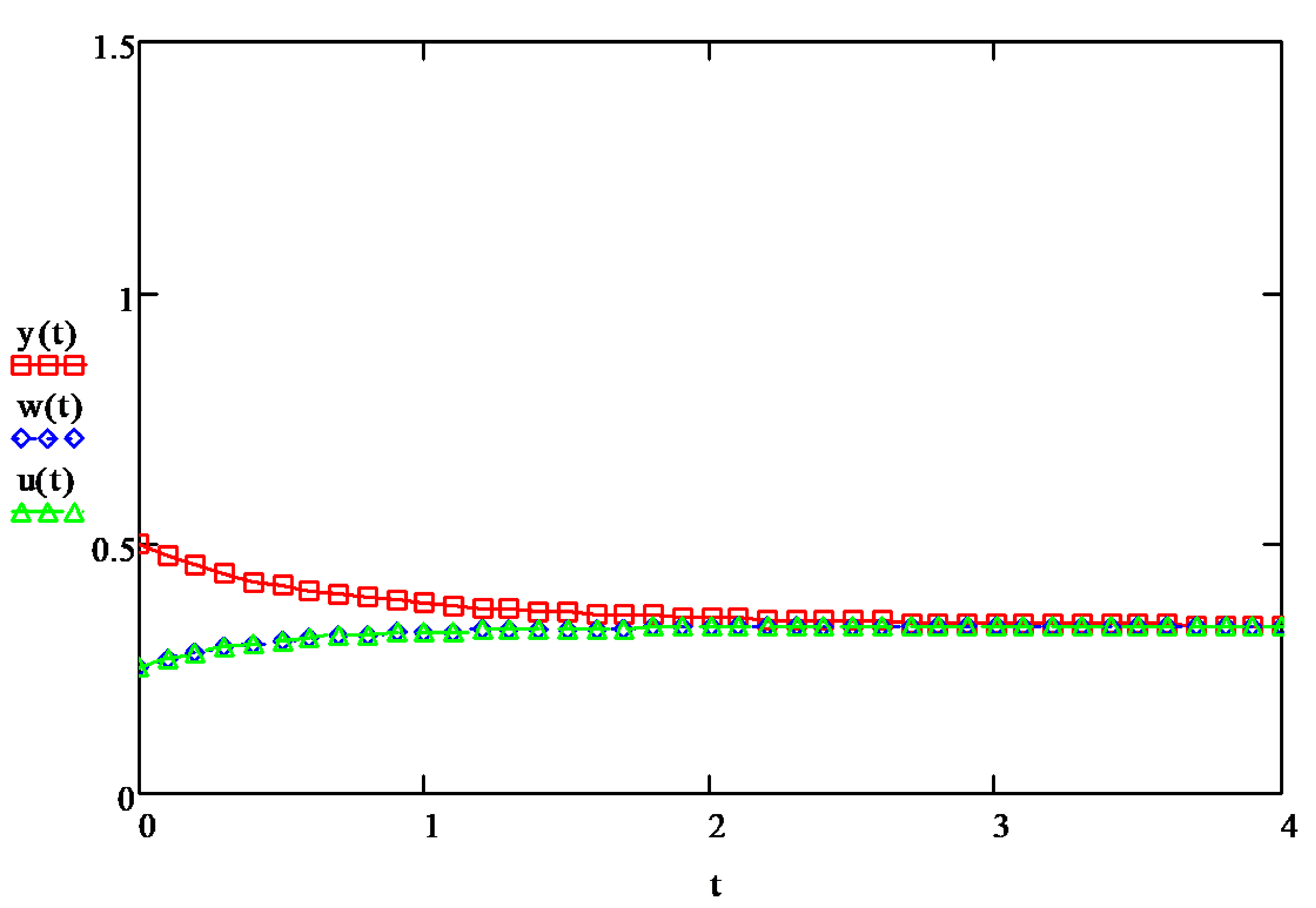}%
\caption{Graph of function \eqref{eq12}, solution of system \eqref{eq10}
with initial values $y(0)=0.5$, $u(0)=0.25$, $w(0)=0.25$.}
\end{center}
\end{figure}

\begin{figure}[H]
\begin{center}
\includegraphics[width=0.5\textwidth]{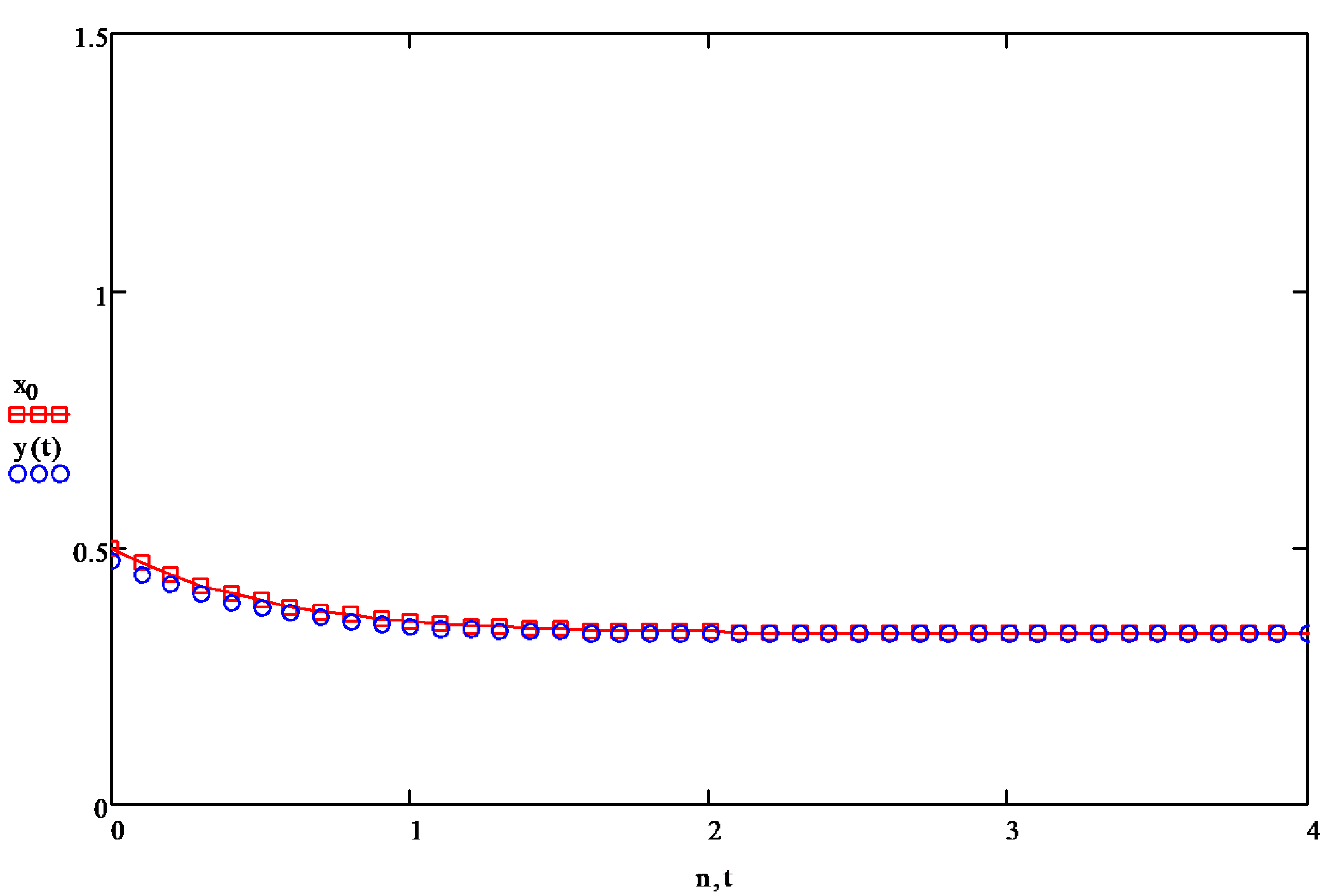}%
\caption{Graph of comparison of the numerical solution of system \eqref{eq1}
$x_{0}(t)$ and analytical solutions \eqref{eq12} $y(t)$ with the same
initial values.}
\end{center}
\end{figure}

\begin{figure}[H]
\begin{center}
\includegraphics[width=0.5\textwidth]{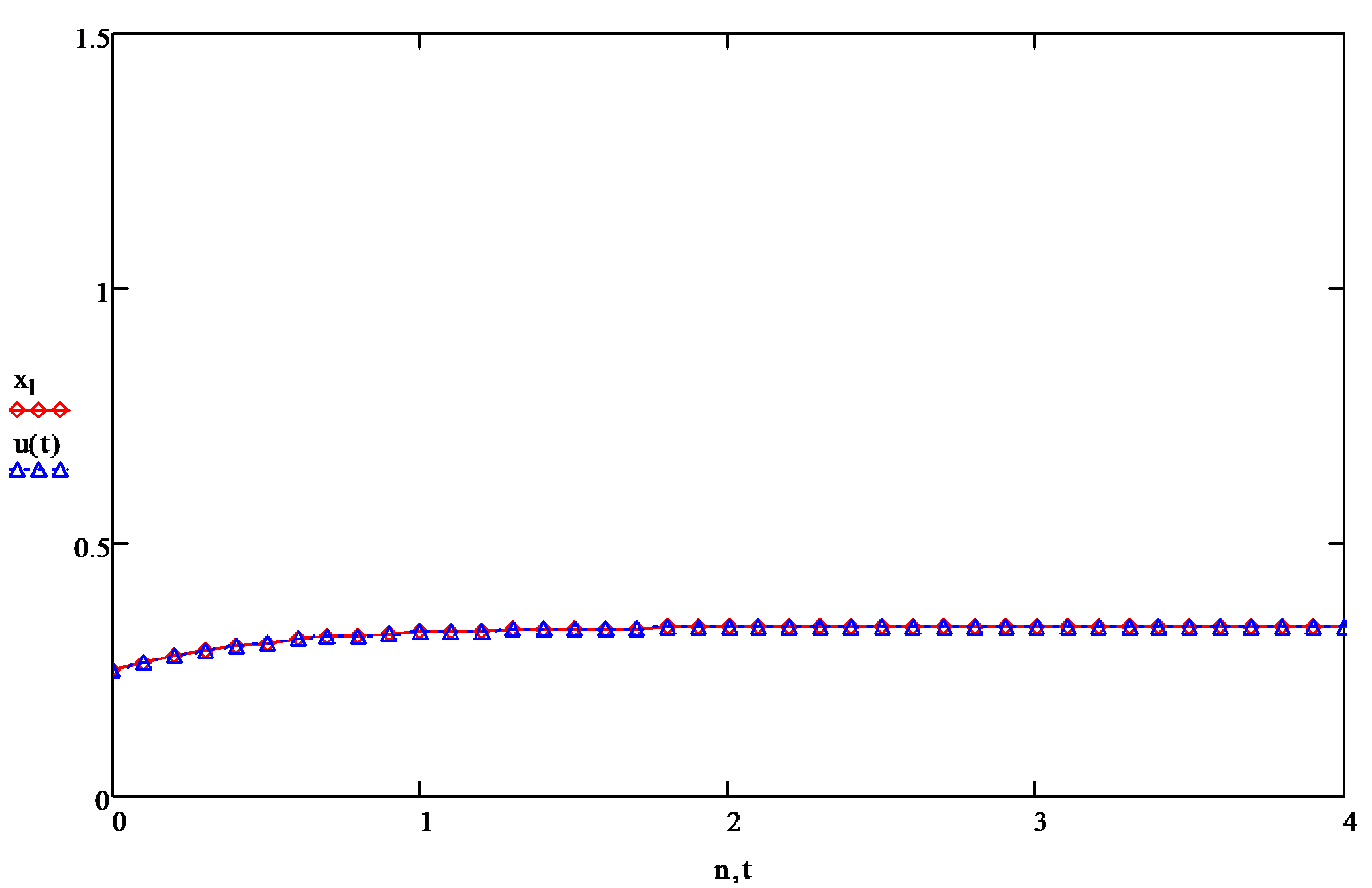}%
\caption{Comparison graph of the numerical solution of system \eqref{eq1} $
x_{1}(t)$ and the analytical solution \eqref{eq12} $w(t)$.}
\end{center}
\end{figure}

A similar picture takes place when comparing the solution of system
\eqref{eq1} with analytical solutions \eqref{eq13} in the line $
\frac{1}{3}<x_{1}(t)$, $x_{2}(t)<1$, $x_{1}(t)=x_{2}(t))$, with the same
initial conditions.

Thereby, the following theorem has been proved.

\textbf{Theorem}. System \eqref{eq1} with the conditions $t\ge 0$, $
x_{0}(t)\ge 0$, $x_{1}(t)\ge 0$, $x_{2}(t)\ge 0$ and $
x_{0}(t)+x_{1}(t)+x_{2}(t)=1$ has a unique fixed point $
(\frac{1}{3},\frac{1}{3},\frac{1}{3})$, which is a hyperbolic saddle.
In addition, the solution \eqref{eq12} and \eqref{eq13} of system \eqref{eq10}
exponentially tends to the solution of the system \eqref{eq1}.

The results of this paper show that the dynamics of an analogue of
strictly non-Volterra operators with continuous time is much richer than
the dynamics of non-Volterra operators with discrete time (see \cite{Mukh15}).
This is also seen from the phase portrait of system \eqref{eq1}, which
consists of five curves. It was found that $x_{0}(t),x_{1}(t),x_{2}(t)
$ tends to the equilibrium point $
M_{2}(\frac{1}{3},\frac{1}{3},\frac{1}{3})$ according to formulas
\eqref{eq12} and \eqref{eq13}, and the sum $x_{0}(t)+x_{1}(t)+x_{2}(t)$ tends to
unity exponentially fast as $t\to +\infty $.

In addition, in the absence of conditions
\[
    x_{0}(t)\ge 0, \quad
    x_{1}(t)\ge 0, \quad
    x_{2}(t)\ge 0
    \quad \textup{and} \quad
    x_{0}(t)+x_{1}(t)+x_{2}(t)=1,
\]
the trajectory defined by system \eqref{eq1} in overall, has several curves.

Each strictly non-Volterra quadratic operator is an interesting example
in the theory of multi\-dimensional nonlinear dynamical systems with
various trajectory behavior.

A comparative analysis of the results obtained in \cite{Mukh15} and in the
present paper shows that the equilibrium positions of system \eqref{eq1}
coincide with the fixed point of the operator \cite{Mukh15} only when $
x_{0}(t)\ge 0$, $x_{1}(t)\ge 0$, $x_{2}(t)\ge 0$ and $
x_{0}(t)+x_{1}(t)+x_{2}(t)=1$ and both of these trajectories tend to
the equilibrium position exponentially fast.

In addition, we can say that considering an analog of a quadratic
operator (studied in \cite{Mukh15}) with continuous time gives some advantage.
Since it has been established that on the basis of computer calculations
it can be said that for $t\ge 4$ solutions of system \eqref{eq1} obtained
using MathCAD (more than 100 solutions with different initial values are
calculated and compared) coincide with the solutions obtained
analytically \eqref{eq12} and \eqref{eq13}.

The author thanks prof. U.A. Rozikov for useful discussions.

%%% References

\EditInfo{December 13, 2020}{March 04, 2021}{Utkir Rozikov}

\end{paper}